\documentclass[letter]{amsart}

\usepackage{amssymb,amsfonts,amsxtra}
\usepackage{paralist}
\usepackage{url}
\usepackage[all]{xypic}
\xyoption{dvips}

\theoremstyle{definition}

\theoremstyle{plain}
\newtheorem{lem}{Lemma}

\newtheorem{thm}{Theorem}
\newtheorem{cor}{Corollary}

\theoremstyle{remark}

\newtheorem{rmk}{Remark}

\newcommand{\CC}{\mathbb{C}}
\newcommand{\GG}{\mathrm{G}}
\newcommand{\OO}{\mathcal{O}}
\newcommand{\PP}{\mathcal{P}}

\newcommand{\TT}{\mathrm{T}}
\newcommand{\ZZ}{\mathbb{Z}}

\renewcommand{\gg}{\mathfrak{g}}
\newcommand{\hh}{\mathfrak{h}}
\newcommand{\mm}{\mathfrak{m}}
\newcommand{\nn}{\mathfrak{n}}

\renewcommand{\tt}{\mathfrak{t}}

\newcommand{\ideal}[1]{{\langle#1\rangle}}

\DeclareMathOperator{\Aut}{Aut}
\DeclareMathOperator{\Der}{Der}
\DeclareMathOperator{\GL}{GL}

\DeclareMathOperator{\Spec}{Spec}
\DeclareMathOperator{\Var}{Var}

\begin{document}

\title[Maximal multihomogeneity of algebraic hypersurface singularities]{Maximal multihomogeneity of\\algebraic hypersurface singularities}

\author{Mathias Schulze}
\address{Mathias Schulze\\
Department of Mathematics\\
Oklahoma State University\\
401 MSCS\\
Stillwater, OK 74078\\
USA}
\email{mschulze@math.okstate.edu}
\thanks{The author is grateful to Michel Granger and Gerhard Pfister for helpful discussions and acknowledges partial support by the Humboldt foundation.}

\date{November 8, 2006}

\begin{abstract}
From the degree zero part of logarithmic vector fields along an algebraic hypersurface singularity we indentify the maximal multihomogeneity of a defining equation in form of a maximal algebraic torus in the embedded automorphism group.
We show that all such maximal tori are conjugate and in one-to-one correspondence to maxmimal tori in the degree zero jet of the embedded automorphism group.

The result is motivated by Kyoji Saito's characterization of quasihomogeneity for isolated hypersurface singularities \cite{Sai71} and extends its formal version \cite[Thm.5.4]{GS06} and a result of Hauser and M\"uller \cite[Thm.4]{HM89}.
\end{abstract}

\subjclass{32S25,17B15}

\keywords{hypersurface singularity, multihomogeneity, logarithmic vector field}

\maketitle

Let $\PP$ be either $\OO=\CC\{x\}$, the ring of convergent power series in $x=x_1,\dots,x_n$, or $\hat\OO=\CC[\![x]\!]$, the ring of formal power series in $x$ and denote by $\mm=\ideal{x}$ the maximal ideal in $\OO$ or $\hat\OO$.
Let $\Aut(\PP)$ be the automorphism group of $\PP$ and let $\Delta(\PP)=\mm\cdot\Der(\PP)=\{\delta\in\Der(\PP)\mid\delta(\mm)\subseteq\mm\}$ be the $\PP$-module of vector fields vanishing at the origin.

Let $0\ne f\in\mm$ and $\Aut_f=\{\varphi\in\Aut(\PP)\mid\varphi(f)\in\ideal{f}\}$ the group of automorphisms preserving the ideal $\ideal{f}$.
In the convergent case, this is the group of automorphisms of $(\CC^n,0)$ preserving the hypersurface $\Var(f)$.
Our object of interest is the $\PP$-module of logarithmic vector fields $\Der_f=\{\delta\in\Der(\PP)\mid\delta(f)\in\ideal{f}\}$ introduced in \cite{Sai80}.
In the convergent case, this is the module of vector fields tangent to the smooth part of $\Var(f)$.
The module $\Der_f$ is unchanged if we assume $f$ to be reduced. 
We shall further assume that $\Der_f\subseteq\Delta(\PP)$.
By Rossi's Theorem \cite[Cor.3.4]{Ros63}, this means in the convergent case that the variety $\Var(f)$ defined by $f$ is not a product with a smooth factor.
As remarked in \cite[2.Rem.(c)]{HM93}, $\Der_f$ is the Lie algebra of the infinite Lie group $\Aut_f$ in the convergent case.

For a fixed coordinate system, any derivation $\delta\in\Der(\PP)$ can be decomposed into homogeneous components, $\delta=\sum_{i=-1}^\infty\delta_i$.
Moreover, $\delta_0=\sum_{i,j}a_{i,j}x_i\partial_{x_j}$ for some matrix $A=(a_{i,j})$ and we call $\delta_0$ diagonal if $A$ is diagonal.
A derivation $\delta\in\Delta(\PP)$ is called nilpotent if $A$ is nilpotent and semisimple if $\mm$ has a basis of eigenvectors of $\delta$.
For a fixed coordinate system, any $\delta\in\Delta(\PP)$ is a sum $\delta=\delta_S+\delta_N$ where $\delta_S=\delta_{S,0}$ and $\delta_{N,0}$ are defined by the semisimple and nilpotent parts of the matrix $A$ corresponding to $\delta_0$.
In particular, $[\delta_S,\delta_{N,0}]=0$.

In \cite[Thm.5.4]{GS06}, K.~Saito's construction \cite[\S3]{Sai71} of the Poincar\'e--Dulac decomposition \cite[Ch.3.\S3.2]{AA88} is generalized to Part \ref{2} of the following theorem.
The statement \ref{2}.\ref{10} is implicitly present in the arguments in \cite{GS06}.
We shall follow the outline of \cite[Thm.4]{HM89} to deduce Part \ref{3} of Theorem \ref{1} from Part \ref{2}.

\begin{thm}[existence of maximal multihomogeneity]\label{1}
Let
\begin{asparaenum}[a.]
\item\label{2} $f\in\hat\OO$ or 
\item\label{3} $f\in\OO$ algebraic over $\CC[x]$
\end{asparaenum}
and let $\delta_1,\dots,\delta_t\in\Der_f$ with diagonal degree $0$ part.
Then there is an algebraic torus $\TT^s\subseteq\Aut_f$ in the sense of \cite[1.Def.ii)]{Mue86} with Lie algebra $\tt^s$ and suitable coordinates such that $\TT^s\subseteq\GL_n(\CC)$.
In these coordinates there is a basis $\sigma_1,\dots,\sigma_s$ of $\tt^s$ which extends to a minimal system of generators $\sigma_1,\dots,\sigma_s,\nu_1,\dots,\nu_r$ of $\Der_f$ and a choice of $f$ with irreducible factors $f_1,\dots,f_m$ such that
\begin{asparaenum}[1.] 
\item\label{12} $\sigma_i$ is diagonal with eigenvalues in $\ZZ$,
\item\label{13} $\nu_i$ is nilpotent,
\item\label{14} $[\sigma_i,\nu_j]\in\ZZ\cdot\nu_j$,
\item\label{15} $\sigma_i(f_j)\in\ZZ\cdot f_j$,
\item\label{10} $(\delta_i)_0\in\langle\sigma_1,\dots,\sigma_s\rangle_\CC$, and
\item\label{11} if $\delta\in\Der_f$ with $[\sigma_i,\delta_0]=0$ for all $i$ then $\delta_S\in\langle\sigma_1,\dots,\sigma_s\rangle_\CC$.
\end{asparaenum}
\end{thm}

\begin{proof}
The statement \ref{2}.\ref{10} holds by the construction in the proof of \cite[Thm.5.4]{GS06} using that the coordinate change in \cite[Thm.5.3]{GS06} is tangent to the identity.

By \ref{2}.\ref{12} and \ref{2}.\ref{15}, there is a formal coordinate change $\bar y(x)$ such that $\bar g(x)=f(\bar y(x))$ fulfils $\sigma_i(\bar g)=\lambda_i\cdot\bar g$ where $\lambda_i\in\ZZ$ and $\sigma_i=\sum_{j=1}^n\sigma_{i,j}x_j\partial_{x_j}$ where $\sigma_{i,j}\in\ZZ$.
We identify $\sigma_i$ with its coefficient vector $(\sigma_{i,1},\dots,\sigma_{i,n})\in\ZZ^n$.
By \cite[II.7.3\,(2)]{Bor91}, the the saturation
\[
L=((\ZZ^n)^\vee/(\ZZ^n/\sum_{i=1}^s\ZZ\sigma_i)^\vee)^\vee\cong\ZZ^s\subseteq\ZZ^n
\]
of the lattice generated by $\sigma_1,\dots,\sigma_s$ defines an algebraic torus
\[
\GG_m^s\cong\TT^s:=\Spec(\CC[L])\subseteq\Spec(\CC[\ZZ^n])\subseteq\GL_n(\CC)\subseteq\Aut(\PP)
\]
with Lie algebra $\tt^s=\langle\sigma_1,\dots,\sigma_s\rangle$.
Since the $x$-monomials are common eigenvectors for $x_1\partial_1,\dots,x_n\partial_n$ with integer eigenvalues, we may assume $L=\sum_{i=1}^s\ZZ\sigma_i$ preserving the condition $\lambda_i\in\ZZ$.
Then the $\sigma$-multihomogeneity of $\bar g$ of multidegree $\lambda$ stated above translates to $\bar g$ being equivariant for $\TT^s\subseteq\GL_n(\CC)$ and the character $\TT^s\to\Spec(\CC[\ZZ])=\GL_1(\CC)$ defined by $L\ni\sigma_i\mapsto\lambda_i\in\ZZ$.
Now $\ideal{f}_{\hat\OO}$ is equivalent to the $\TT^s$-stable ideal $\ideal{\bar g}_{\hat\OO}$ and \cite[Thm.2']{HM89} implies that also $\ideal{f}_\OO$ is equivalent to a $\TT^s$-stable ideal.
This means that there is an analytic coordinate change $y(x)$ and a unit $u\in\OO^*$ such that $g(x)=u(x)\cdot f(y(x))$ generates a $\TT^s$-stable ideal in $\OO$.
By abuse of notation we denote this $g$ by $f$ again.
Then \cite[Hilfssatz~2]{Mue86} shows that $\ideal{f}_\OO$ has a $\TT^s$-equivariant generator which we may assume to be $f$.
This proves \ref{3}.\ref{12} and $\sigma_i(f)\in\ZZ\cdot f$ which is weaker than \ref{3}.\ref{15}.

To avoid the non-trivial \cite[Hilfssatz~2]{Mue86} and conclude \ref{3}.\ref{15}, one can argue as follows:
By \cite[Thm.13.2]{Hum75}, the $\TT^s$- and $\tt^s$-stable subspaces in $\OO/\mm^k$ coincide.
Since ideals in $\OO$ are closed in the $\mm$-adic topology, this shows that $\ideal{f}_\OO$ is $\sigma$-stable.
Then a multigraded version\footnote{The arguments in \cite[2.2-4]{SW73} give a more general correspondence of $(k^s,+)$-graduations and sets of $s$ simultaneously diagonalizable $k$-derivations on analytic $k$-algebras.} of \cite[2.4]{SW73} shows that $\ideal{f}_\OO$ has a $\TT^s$-equivariant generator.
With $\ideal{f}_\OO$ also its minimal associated primes are $\sigma$-stable by \cite[2.5]{SW73}.
Again by \cite[2.4]{SW73}, these primes have $\TT^s$-equivariant generators $f_1,\dots,f_m$ and \ref{3}.\ref{15} follows.

With $f$ also its partial derivatives $\partial_{x_i}(f)$ are $\TT^s$-equivariant.
Thus $\Der_f$ is $\TT^s$-stable since it can be considered as the projection of the syzygy module of $\partial_{x_1}(f),\dots,\partial_{x_n}(f),f$ to the first $n$ components.
Its $\TT^s$-sub-module $\tt^s=\langle\sigma_1,\dots,\sigma_s\rangle$ is the adjoint representation of $\TT^s$ and hence trivial.
By a module version\footnote{The statement in \cite[Hilfssatz~2]{Mue86} holds more generally for any analytic sub-module of a free analytic module in the sense of \cite{GR71}.} of \cite[Hilfssatz~2]{Mue86}, $\Der_f$ has a minimal system of generators $\sigma_1,\dots,\sigma_s,\nu_1,\dots,\nu_r$ which spans over $\CC$ a rational $\TT^s$-module.
By the proof of \cite[Thm.9.13]{Mil06} the $\TT^s$-module $\langle\sigma_1,\dots,\sigma_s,\nu_1,\dots,\nu_r\rangle_\CC$ can be diagonalized without changing the $\sigma_i$ and \ref{3}.\ref{14} follows.
Like above, one can use a multigraded module version of \cite[2.4]{SW73} instead.

The maximality property \ref{3}.\ref{11} follows from its formal version \ref{2}.\ref{11}.
Combined with \ref{3}.\ref{14} it guarantees the existence of $\nu_i$ satisfying \ref{3}.\ref{13}:
If the $\sigma$-multidegree of a $\nu_i$ is non-zero then it is nilpotent by \cite[Lem.2.6]{GS06}, otherwise one can subtract its semisimple part which is a linear combination of the $\sigma_i$ by \ref{3}.\ref{11}.
\end{proof}

\begin{rmk}\
\begin{asparaenum}
\item By \cite[2.11]{KPR75}, the Implicit Function Theorem holds in $\CC\langle x\rangle$, the ring of algebraic power series.
Thus the proof of \cite[Thm.2']{HM89} works inside of $\CC\langle x\rangle$ as well.
However it is not clear if this ring is graded in the sense of \cite[\S1-2]{SW73}.
A positive answer would imply that Theorem \ref{1}.\ref{3} even holds inside of $\CC\langle x\rangle$.
\item According to Michel Granger, the uniqueness of $s$ in Theorem \ref{1} can be deduced from the conjugacy of all Cartan subalgebras \cite[III.4.Thm.2]{Ser87} by showing that $\tt^s$ and $\nn_0=\ideal{(\nu_j)_0\mid\sigma_i(\nu_j)=0}_\CC$ span a Cartan subalgebra $\hh$ of $\gg=\Der_f/(\Der_f\cap\mm^2\cdot\Der(\PP))$ where $\nn_0$ is the intersection of $\hh$ with the set of nilpotent elements in $\gg$.
\end{asparaenum}
\end{rmk}

We shall use results of M\"uller \cite{Mue86} to prove a stronger statement than uniqueness of $s$.
Consider the group morphisms $\pi_k:\Aut(\PP)\to\Aut_k(\PP)=\Aut(\PP/\mm^{k+2})$ and the Lie algebra morphisms $\pi_k:\Delta(\PP)\to\Delta_k(\PP)=\Delta(\PP)/(\mm^{k+1}\cdot\Delta(\PP))$.
Note that $\Aut_k(\PP)$ is an algebraic group with Lie algebra $\Delta(\PP)$.
Like in \cite[\S2]{Mue86}, one can use Artin's Approximation Theorem \cite[Thm.1.2]{Art68} to prove the first part of the following

\begin{lem}\label{4}
$\pi_k(\Aut_f)$ is an algebraic group with Lie algebra $\pi_k(\Der_f)$.
\end{lem}

\begin{proof}
By exactness of completion $\Der_{\hat f}=\widehat\Der_f$ where $\hat f$ denotes $f$ considered in $\hat\OO$.
Thus $\pi_k(\Der_f)=\pi_k(\Der_{\hat f})$ and we may assume that $\PP=\hat\OO$.
Consider the Lie algebra morphisms $\pi_k^m:\Delta_m(\PP)\to\Delta_k(\PP)$ and denote by $f_k$ the image of $f$ in $\OO/\mm^{k+1}$.
For fixed $k$, the $\pi_k^m(\Der_{f_m})$ form a decreasing sequence of subvector spaces in $\Delta_k(\PP)$ which implies $\pi_k^m(\Der_{f_m})=D_k$ for large $m$.
Then $\pi_m^{m+1}:D_{m+1}\to D_m$ is surjective and hence $D_k=\pi_k(\Der_f)$.
In the proof of the first statement, we find $\pi_k^m(\Aut_{f_m})=\pi_k(\Aut_f)$ for large $m$ and $\Der_{f_m}$ is the Lie algebra of $\Aut_{f_m}$.
\end{proof}

\begin{thm}[uniqueness of maximal multihomogeneity]\label{5}
The algebraic torus $\TT^s$ in Theorem \ref{1} is maximal in $\Aut_f$ and also $\pi_0(\TT^s)\subseteq\pi_0(\Aut_f)$ is a maximal algebraic torus.
All maximal algebraic tori in $\Aut_f$ and $\pi_0(\Aut_f)$ are conjugate.
\end{thm}

\begin{proof}
First let $\TT^t\subseteq\pi_0(\Aut_f)$ be an algebraic torus with Lie algebra $\tt^t$.
If $\pi_0(\TT^s)\subseteq\TT^t$ then the Lie algebra of $\TT^t$ consists of $\sigma$-homogeneous semisimple elements of multidegree $0$ and hence $\pi_0(\tt^s)=\tt^t$ by Theorem \ref{1}.\ref{11}.
By Lemma \ref{4} and \cite[Cor.13.11]{Mil06}, this implies $\pi_0(\TT^s)=\TT^t$ and hence $\pi_0(\TT^s)$ is maximal.

Now let $\TT^t\subseteq\Aut_f$ be an algebraic torus in the sense of \cite[1.Def.ii)]{Mue86} with Lie algebra $\tt^t$.
By definition, $\pi_k(\TT^t)\subseteq\pi_k(\Aut_f)\subseteq\Aut_k(\PP)$ defines a rational representation of $\TT^t$ which is diagonalizable by \cite[Thm.9.13]{Mil06}. 
Then Cartan's Uniqueness Theorem holds for $\TT^t$ by \cite[Satz]{Kau67} and $\pi_0$ restricted to $\TT^t$ is injective.
This shows that maximality of $\pi_0(\TT^s)$ implies maximality of $\TT^s$.

By \cite[VIII.21.3.Cor.A]{Hum75}, all maximal tori of $\pi_0(\Aut_f)$ are conjugate.
An embedded version\footnote{In the first paragraph of \cite[\S5]{Mue86}, the proof of \cite[Satz 2]{Mue86} is reduced to the case of embedded automorphism groups using \cite[Satz 6]{Mue86}.} of \cite[Satz 2]{Mue86} shows that a conjugacy in $\pi_0(\Aut_f)$ of algebraic tori in $\Aut_f$ lifts to a conjugacy in $\Aut_f$.
\end{proof}

\begin{cor}[lifting of maximal tori]
$\pi_0:\Aut_f\to\pi_0(\Aut_f)$ defines a bijection of algebraic tori.
\end{cor}

\bibliographystyle{amsalpha}
\bibliography{mmahs}

\end{document}